\documentclass[a4paper,11pt]{article}

\def\qed{\hfill $\squar$}
\def\squar{\vbox{\hrule\hbox{\vrule height 6pt \hskip 6pt\vrule}\hrule}}

\usepackage{amsmath}
\usepackage{amssymb}
\usepackage{epsf}
\usepackage{graphicx,amssymb,amsmath,color,latexsym} 

\newcommand{\eps}{\varepsilon}

\newfont{\gothic}{eufm10}

\newtheorem{theorem}{Theorem}[section]
\newtheorem{lemma}[theorem]{Lemma}
\newtheorem{proposition}[theorem]{Proposition}

\newtheorem{definition1}[theorem]{Definition}

\newtheorem{example1}[theorem]{Example}

\def\barray{\begin{eqnarray*}}             \def\earray{\end{eqnarray*}}
\def\beq{\begin{equation}} \def\eeq{\end{equation}}

\usepackage{times}

\begin{document}
\renewcommand{\figurename}{\footnotesize Fig.}

\begin{center}
{\LARGE  Bifurcation of asymptotically stable periodic solutions in nearly impact  oscillators: DRAFT }\\
\bigskip

{\large  Oleg Makarenkov}\\

\vskip0.2truecm
Department of Mathematics, Imperial College London, UK \\

\

{\large  Ferdinand Verhulst}\\

\vskip0.2truecm Mathematical Institute, University of Utrecht,
Netherlands

\

09 September 2009

\end{center}

\section{Introduction}

 Since the observation by Glover-Lazer-McKenna
 \cite{lazer} that a simple harmonic oscillator with a piecewise
 linear stiffness (jumping nonlinearity) contributes to the explanation of

 the failure
 of the Takoma bridge, the studying of periodic oscillations in such
 models got a lot of attention of mathematicians; see the recent
 survey \cite{mawhin} and the papers \cite{yagasaki,blm}. Also, new
 engineering studies of impact oscillators open up a large potential
 for challenging extensions of these results. In fact, Ivanov
  \cite{ivanov} argued that harmonic oscillators with a
 jumping nonlinearity with one part of the force field nearly infinite is a better
 model for describing the bouncing ball, rather then its limit
 version for an impact oscillator. In our modeling the resulting system of
 differential equations is singularly perturbed, but as we discuss
 below, the classical singular perturbation theory does not apply.
 In this paper we develop an averaging-like approach which solves
 the problem in a weakly nonlinear case. For a discussion of the
 use of averaging method for regular impacting systems we refer to
 \cite{BK01}. To be explicit, our approach concerns the existence and
 stability of periodic oscillations of the following system
 \begin{equation}\label{ps}
 \begin{array}{l}
 \ddot{x}+x=\eps f(t,x,\dot{x},\eps), \qquad x>0, \\
 \ddot{x}+\dfrac{1}{\eps^2(\omega_\eps)^2}x=g(t,x,\dot{x},\eps),
 \quad x<0,
 \end{array}
 \end{equation}
 where $f,g \in C^1(\mathbb{R} \times \mathbb{R}\times
 \mathbb{R}\times[0,1],\mathbb{R}),\eps>0$ is a small parameter,
 $\omega_\eps \to \omega_0 \in \mathbb{R} \ {\rm as}\ \eps \to 0.$
 System \ref{ps} can be considered as a smoothed version of a system with impacts.
 We will study  resonance oscillations and assume,
 therefore, that
 \begin{equation*}
 f(t+\pi,u,v,\eps) \equiv f(t,u,v, \eps), g(t+\pi,u,v,\eps) \equiv g(t,u,v,\eps).
 \end{equation*}
 System (\ref{ps}) represents a natural singular perturbation description of impact phenomena which is different from the usual approaches.
 Our main result (Theorem 1) states that the emergence of asymptotically stable $\pi$-periodic
 solutions in (\ref{ps}) from $\pi$-periodic cycles of non-smooth limiting the system
 \begin{equation}\label{np}
 \begin{array}{l}
 \ddot{x}+x=0, \qquad x>0,\\
 \dot{x}(t-0)=-\dot{x}(t+0), \quad x(t)=0,
 \end{array}
 \end{equation}
 can be studied by a special form of the averaging method combined
 with a suitable scaling of time when solutions pass the half plane
 $x<0$. This involves the use of the implicit function theorem for
 a non-smooth problem in the limit as $\varepsilon \rightarrow 0$;
 this is possible by introducing a suitable Poincar\'e map. The
 result is a change of stability of a fixed point when its
 eigenvalues enter the  unit disc from
 outside through the imaginary axis.

\section{Main result}
We prove the following theorem.
\begin{theorem}\nonumber
Let $f,g \in C^1(\mathbb{R}\times\mathbb{R}\times\mathbb{R},\mathbb{R})$ be $\pi$-periodic with respect to time. Define
\begin{displaymath}
\begin{array}{rcl}
\hskip-0.2cm\overline{P}(A,\theta)&=&-\int\limits_0^{\frac{\pi}{2}-\theta}\begin{pmatrix} \sin(\tau+\theta)\\ \dfrac{1}{A}\cos(\tau+\theta)\end{pmatrix} (f(\tau,A\cos(\tau+\theta),-A\sin(\tau+\theta),0)-2\omega_0A\cos(\tau+\theta))d\tau-\\
&&\hskip-0.3cm-\int\limits_{\frac{\pi}{2}-\theta}^\pi \begin{pmatrix} \sin(\tau+\theta+\pi)\\ \dfrac{1}{A}\cos(\tau+\theta+\pi)\end{pmatrix} (f(\tau,A\cos(\tau+\theta+\pi),-A\sin(\tau+\theta+\pi),0)-\\
&&\hskip-0.3cm-2\omega_0A\cos(\tau+\theta+\pi))d\tau-\omega_0
\int\limits_0^\pi \begin{pmatrix}
\sin\left(s+\dfrac{\pi}{2}\right)\\0\end{pmatrix}g \left(
\dfrac{\pi}{2}-\theta,0,-A\sin \left(
s+\dfrac{\pi}{2}\right),0\right)ds.
\end{array}
\end{displaymath}
If $\overline{P}(A_0,\theta_0)=0$ for some $(A_0,\theta_0)\in\mathbb{R}\times(0,\pi)$ and the real parts of eigenvalues of $\overline{P}'(A_0,\theta_0)$ are negative, then, for any $\eps>0$ sufficiently small, equation (\ref{ps}) has a unique $\pi$-periodic solution satisfying
\begin{equation*}
x_\eps(t) \to x_0(t)\quad \mbox{as}\quad \eps \to 0\quad
\mbox{pointwise on} \quad [0,\pi]\setminus \left\{
\dfrac{\pi}{2}-\theta_0 \right\}
\end{equation*}
where $x_0$ is the unique $\pi$-periodic solution of the equation (\ref{np}) with initial condition
\begin{equation*}
(x(0),\dot{x}(0))=(A_0\cos \theta_0,-A_0\sin \theta_0).
\end{equation*}
Moreover, the solution $x_\eps$ is asymptotically stable.
\end{theorem}

\begin{figure}\label{fig1}
\begin{center}
\includegraphics[scale=0.9]{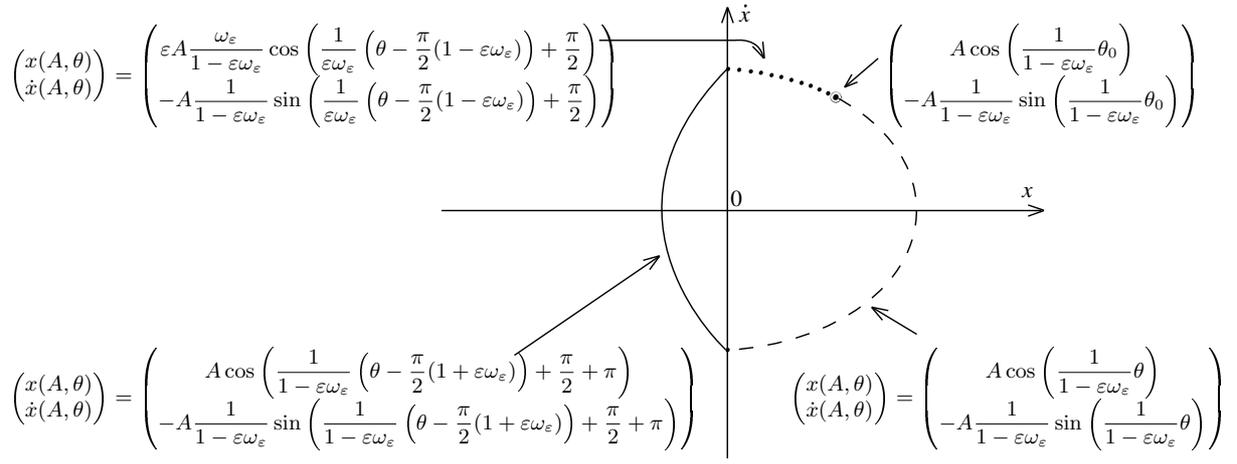}
\footnotesize\caption{Illustration of the change of variables
(\ref{CH1})-(\ref{CH3}).}
\end{center}
\end{figure}

\noindent{\bf Proof.} Rewrite system (\ref{ps}) as follows (see
Fig.~\ref{fig1})
\begin{eqnarray}
&&\ddot{x}+\dfrac{1}{(1-\eps\omega_\eps)^2}x=\eps f(t,x,\dot{x},\eps)-2\eps\dfrac{\omega_\eps}{(1-\eps\omega_\eps)^2}x-\eps^2\dfrac{(\omega_\eps)^2}{(1-\eps\omega_\eps)^2}x, \quad x>0, \label{EQ1}\\
&&\ddot{x}+\dfrac{1}{(\eps)^2(\omega_\eps)^2}x=g(t,x,\dot{x},\eps),
\quad x<0,\label{EQ2}
\end{eqnarray}
so, that any solution of the reduced system ($\varepsilon =0$)
\begin{eqnarray*}
&&\ddot{x}+\dfrac{1}{(1-\eps\omega_\eps)^2}x=0, \quad x>0,\\
&&\ddot{x}+\dfrac{1}{(\eps)^2(\omega_\eps)^2}x=0, \quad x<0
\end{eqnarray*}
is $\pi$-periodic. Let us introduce variables $(A,\theta)$ as follows
\begin{equation}\label{CH1}
\left\{
\begin{array}{rcll}
x&=&A\cos\left( \dfrac{1}{1-\eps\omega_\eps}\theta\right),&\\
\dot{x}&=&-A \dfrac{1}{1-\eps\omega_\eps}\sin \left(
\dfrac{1}{1-\eps\omega_\eps}\theta\right), & \theta \in
\left[0,\dfrac{\pi}{2}(1-\eps\omega_\eps)\right],
\end{array}
\right.
\end{equation}
\begin{equation}\label{CH2}
\hskip-0.15cm\left\{
\begin{array}{rcll}
x&=&\eps A\dfrac{\omega_\eps}{1-\eps\omega_\eps}\cos\left( \dfrac{1}{\eps\omega_\eps}(\theta-\dfrac{\pi}{2}(1-\eps\omega_\eps))+\dfrac{\pi}{2}\right),&\\
\dot{x}&=&-A \dfrac{1}{1-\eps\omega_\eps}\sin \left(
\dfrac{1}{\eps\omega_\eps}(\theta-\dfrac{\pi}{2}(1-\eps\omega_\eps))+\dfrac{\pi}{2}\right),
& \theta \in
\left[\dfrac{\pi}{2}(1-\eps\omega_\eps),\dfrac{\pi}{2}(1+\eps\omega_\eps)\right],
\end{array}
\right.
\end{equation}
\begin{equation}\label{CH3}
\hskip-0.15cm\left\{
\begin{array}{rcll}
x&=&A\cos \left(\dfrac{1}{1-\eps\omega_\eps}\left( \theta -\dfrac{\pi}{2}(1+\eps\omega_\eps)\right)+\dfrac{\pi}{2}+\pi\right),&\\
\dot{x}&=&-A \dfrac{1}{1-\eps\omega_\eps}\sin \left(
\dfrac{1}{1-\eps\omega_\eps}\left(\theta-\dfrac{\pi}{2}(1+\eps\omega_\eps)\right)+\dfrac{\pi}{2}+\pi\right),
& \theta \in \left[\dfrac{\pi}{2}(1+\eps\omega_\eps),\pi\right],
\end{array}
\right.
\end{equation}
which transforms equations (\ref{EQ1})-(\ref{EQ2}) to the following system
\begin{eqnarray}
&&\begin{pmatrix} \dot{A}\\ \dot{\theta}\end{pmatrix}=\begin{pmatrix} 0\\1\end{pmatrix}+\eps G_1(t,A,\theta,\eps), \qquad \mbox{if} \quad \theta \in \left[0,\dfrac{\pi}{2}(1-\eps\omega_\eps)\right], \label{S1}\\
&&\begin{pmatrix} \dot{A}\\ \dot{\theta}\end{pmatrix}=\begin{pmatrix} 0\\1\end{pmatrix}+\eps G_2(t,A,\theta,\eps), \qquad \mbox{if} \quad \theta \in \left[\dfrac{\pi}{2}(1-\eps\omega_\eps),\dfrac{\pi}{2}(1+\eps\omega_\eps)\right], \label{S2}\\
&&\begin{pmatrix} \dot{A}\\
\dot{\theta}\end{pmatrix}=\begin{pmatrix} 0\\1\end{pmatrix}+\eps
G_3(t,A,\theta,\eps), \qquad \mbox{if} \quad \theta \in
\left[\dfrac{\pi}{2}(1+\eps\omega_\eps),\pi\right], \label{S3}
\end{eqnarray}
where
\begin{equation*}
G_1(t,A,\theta,\eps)= \left(
\begin{array}{l}
-(1-\eps\omega_\eps)\sin \left( \dfrac{1}{1-\eps\omega_\eps}\theta \right)\left( f \left( t,A \cos \left( \dfrac{\theta}{1-\eps\omega_\eps} \right),\right.\right.\\
\left.\quad \qquad\qquad \qquad \qquad\qquad \qquad-A \dfrac{1}{1-\eps\omega_\eps}\sin \left( \dfrac{\theta}{1+\eps\omega_\eps} \right),\eps \right)-\\
\quad \qquad\qquad \qquad \qquad\qquad \qquad \left. -A \dfrac{\omega_\eps}{(1-\eps\omega_\eps)^2}(2+\eps\omega_\eps)\cos \left( \dfrac{\theta}{1-\eps\omega_\eps} \right) \right)\\
-\dfrac{(1-\eps\omega_\eps)^2}{A}\cos \left( \dfrac{1}{1-\eps\omega_\eps}\theta \right)\left( f \left( t,A \cos \left( \dfrac{\theta}{1-\eps\omega_\eps} \right),\right.\right.\\
\qquad \qquad\qquad \qquad\qquad \qquad\quad\left.-A \dfrac{1}{1-\eps\omega_\eps}\sin \left( \dfrac{\theta}{1+\eps\omega_\eps} \right),\eps \right)-\\
\qquad \qquad\qquad \qquad\qquad \qquad\quad\left. -A
\dfrac{\omega_\eps}{(1-\eps\omega_\eps)^2}(2+\eps\omega_\eps)\cos
\left( \dfrac{\theta}{1-\eps\omega_\eps} \right) \right)
\end{array}\right)
\end{equation*}
\begin{equation*}
G_2(t,A,\theta,\eps)= \left(
\begin{array}{l}
-\dfrac{1}{\eps}(1-\eps\omega_\eps)\sin \left( \dfrac{1}{\eps\omega_\eps}\left( \theta -\dfrac{\pi}{2}(1-\eps\omega_\eps)\right)+\dfrac{\pi}{2} \right) \cdot\\
\qquad \ \cdot g \left( t,\eps A \dfrac{\omega_\eps}{1-\eps\omega_\eps}\cos \left(\dfrac{1}{\eps\omega_\eps}\left( \theta-\dfrac{\pi}{2}(1-\eps\omega_\eps) \right) +\dfrac{\pi}{2}\right),\right.\\
\qquad \qquad \left.-A \dfrac{1}{1-\eps\omega_\eps}\sin \left( \dfrac{1}{\eps\omega_\eps} \left( \theta-\dfrac{\pi}{2}(1-\eps\omega_\eps) \right)+\dfrac{\pi}{2}\right),\eps \right)\\
-\dfrac{1}{A}(1-\eps\omega_\eps)\omega_\eps \cos \left( \dfrac{1}{\eps\omega_\eps} \left( \theta-\dfrac{\pi}{2}(1-\eps\omega_\eps) \right) +\dfrac{\pi}{2}\right) \cdot\\
\qquad \ \cdot g \left( t,\eps A \dfrac{\omega_\eps}{1-\eps\omega_\eps}\cos \left(\dfrac{1}{\eps\omega_\eps}\left( \theta-\dfrac{\pi}{2}(1-\eps\omega_\eps) \right) +\dfrac{\pi}{2}\right),\right.\\
\qquad \qquad \left.-A \dfrac{1}{1-\eps\omega_\eps}\sin \left(
\dfrac{1}{\eps\omega_\eps} \left(
\theta-\dfrac{\pi}{2}(1-\eps\omega_\eps)
\right)+\dfrac{\pi}{2}\right),\eps \right)
\end{array}\right)
\end{equation*}
\begin{equation*}
G_3(t,A,\theta,\eps)= \left(
\begin{array}{l}
-(1-\eps\omega_\eps)\sin \left( \dfrac{1}{1-\eps\omega_\eps}\left( \theta -\dfrac{\pi}{2}(1+\eps\omega_\eps)\right)+\dfrac{\pi}{2} +\pi\right) \cdot\\
\quad \cdot \left( f \left( t,A \cos \left(\dfrac{1}{1-\eps\omega_\eps}\left( \theta-\dfrac{\pi}{2}(1+\eps\omega_\eps) \right) \right)+\dfrac{\pi}{2}+\pi\right),\right.\\
\quad \quad \left.-A \dfrac{1}{1-\eps\omega_\eps}\sin \left( \dfrac{1}{1-\eps\omega_\eps} \left( \theta-\dfrac{\pi}{2}(1+\eps\omega_\eps) \right)+\dfrac{\pi}{2}+\pi \right),\eps \right)-\\
\quad\quad\quad-A\dfrac{\omega_\eps}{(1-\eps\omega_\eps)^2}(2+\eps\omega_\eps)\cos\left( \dfrac{1}{1-\eps\omega_\eps} \left( \theta-\dfrac{\pi}{2}(1+\eps\omega_\eps) \right) +\dfrac{\pi}{2}+\pi \right)\\
-\dfrac{(1-\eps\omega_\eps)^2}{A}\cos \left( \dfrac{1}{1-\eps\omega_\eps} \left( \theta-\dfrac{\pi}{2}(1-\eps\omega_\eps) \right) +\dfrac{\pi}{2}+\pi\right) \cdot\\
\quad\cdot \left( f \left( t,A \cos \left(\dfrac{1}{1-\eps\omega_\eps}\left( \theta-\dfrac{\pi}{2}(1+\eps\omega_\eps) \right) \right)+\dfrac{\pi}{2}+\pi\right),\right.\\
\quad \quad \left.-A \dfrac{1}{1-\eps\omega_\eps}\sin \left( \dfrac{1}{1-\eps\omega_\eps} \left( \theta-\dfrac{\pi}{2}(1+\eps\omega_\eps) \right)+\dfrac{\pi}{2}+\pi \right),\eps \right)-\\
\quad\quad\quad-A\dfrac{\omega_\eps}{(1-\eps\omega_\eps)^2}(2+\eps\omega_\eps)\cos\left(
\dfrac{1}{1-\eps\omega_\eps} \left(
\theta-\dfrac{\pi}{2}(1+\eps\omega_\eps) \right)
+\dfrac{\pi}{2}+\pi \right)\end{array}\right).
\end{equation*}
If $t \mapsto (A(t),\theta(t)-t)$ is an asymptotically stable
$\pi$-periodic solution of (\ref{S1})-(\ref{S3}) then
$(x,\dot{x})$ defined by (\ref{CH1})-(\ref{CH3}) is an
asymptotically stable $\pi$-periodic solution of
(\ref{EQ1})-(\ref{EQ2}). To prove the existence of asymptotically
stable $\pi$-periodic solutions of equations (\ref{S1})-(\ref{S3}) we show
that each solution\linebreak
$(\overline{A}(\cdot,A,\theta,\eps),\overline{\theta}(\cdot,A,\theta,\eps))$
of (\ref{S1})-(\ref{S3}) with initial condition
$(\overline{A}(0,A,\theta,\eps),\overline{\theta}(0,A,\theta,\eps))=(A,\theta)$
is defined on $[0,\pi]$ whenever $(A,\theta)$
belongs to a small neighborhood of $(A_0,\theta_0)$, and that the
map
\begin{equation}\label{PO}
P_\eps(A,\theta)=(\overline{A}(\pi,A,\theta,\eps),\overline{\theta}(\pi,A,\theta,\eps)-\pi)
\end{equation}
contracts in this neighborhood.

{\bf Step 1.} First we show that solution $t \mapsto (\overline{A}(\cdot,A,\theta,\eps),\overline{\theta}(\cdot,A,\theta,\eps))$ of (\ref{S1})-(\ref{S3}) on $[0,\pi]$ can be consequently sewed by solutions of systems (\ref{S1}), (\ref{S2}) and (\ref{S3}).

Denote by $t \mapsto
(\overline{A}_i(\cdot,t_0,A,\theta,\eps),\overline{\theta}_i(\cdot,t_0,A,\theta,\eps)),$
$i=1,2,3,$ the solutions of (\ref{S1}), (\ref{S2}), (\ref{S3})
respectively with initial condition
$(\overline{A}_i(t_0,t_0,A,\theta,\eps),\overline{\theta}_i(t_0,t_0,A,\theta,\eps))=(A,\theta).$
Put
\begin{equation*}
F_1(T,A,\theta,\eps)=\dfrac{1}{1-\eps\omega_\eps}\overline{\theta}_1(T,0,A,\theta,\eps)-\dfrac{\pi}{2}.
\end{equation*}
Since
$F_1\left(\dfrac{\pi}{2}-\theta_0,A_0,\theta_0,0\right)=\overline{\theta}_1\left(\dfrac{\pi}{2}-\theta_0,0,A_0,\theta_0,0\right)-\dfrac{\pi}{2}=0$
and
$$(F_1)'_T\left(\dfrac{\pi}{2}-\theta_0,A_0,\theta_0,0\right)=(\overline{\theta})'_{(1)}\left(\dfrac{\pi}{2}-\theta_0,0,A_0,\theta_0,0\right)=1$$
then by the implicit function theorem \cite[Ch.~X, \S~2,
Theorems~1 and 2]{kolm} there exists $T_1 \in
C^1(\mathbb{R}\times\mathbb{R}\times [0,1],\mathbb{R})$ such that
\begin{equation*}
F_1(T_1(A,\theta,\eps),A,\theta,\eps)=0, \qquad |A-A_0|<\delta,\
|\theta-\theta_0|<\delta, \ \eps \in [0,\delta),
\end{equation*}
where $\delta>0$ sufficiently small. Or, equivalently,
\begin{equation*}
\dfrac{1}{1-\eps\omega_\eps}\overline{\theta}_1(T_1(A,\theta,\eps),0,A,\theta,\eps)=\dfrac{\pi}{2},
\qquad |A-A_0|<\delta,\ |\theta-\theta_0|<\delta,\ \eps \in
[0,\delta).
\end{equation*}
Therefore, the solution of system (\ref{S1}) with initial
condition $(A,\theta)$ at $t=0$ approaches the threshold of
switching to (\ref{S2}) at time $T_1(A,\theta,\eps).$

Now we show that the solution
$$
\begin{array}{l}
\left(\overline{A}_2\left(\cdot,T_1(A,\theta,\eps),\overline{A}_1(T_1(A,\theta,\eps),A,\theta,\eps),
\dfrac{\pi}{2}(1-\eps\omega_\eps),\eps\right),\right.\\
\qquad
\left.\overline{\theta}_2\left(\cdot,T_1(A,\theta,\eps),\overline{A}_1(T_1(A,\theta,\eps),A,\theta,\eps),\dfrac{\pi}{2}(1-\eps\omega_\eps),\eps\right)\right)
\end{array}
$$
stays till some time $T_2(A,\theta,\eps)$ in
\begin{equation*}
[0,\infty)\times\left[\dfrac{\pi}{2}(1-\eps\omega_\eps),\dfrac{\pi}{2}(1+\eps\omega_\eps)\right]
\end{equation*}
and that $T_2(A,\theta,\eps)$ is given by
\begin{equation*}
T_2(A,\theta,\eps)=T_1(A,\theta,\eps)+\eps\widetilde{T}_2(A,\theta,\eps),
\end{equation*}
where $\widetilde{T}_2\in C^1(\mathbb{R}\times\mathbb{R}\times[0,1],\mathbb{R})$. To do this consider
\begin{eqnarray*}
F_2(T,A,\theta,\eps)&=&\dfrac{1}{\eps\omega_\eps} \left( \overline{\theta}_2 (T_1(A,\theta,\eps)+\eps T,T_1(A,\theta,\eps),  \right.\\
&&\left.\overline{A}_1(T_1(A,\theta,\eps),0,A,\theta,\eps),\overline{\theta}_1(T_1(A,\theta,\eps),0,A,\theta,\eps),\eps)-\dfrac{\pi}{2}(1-\eps\omega_\eps)\right)-\pi.
\end{eqnarray*}
\begin{equation*}
F_2(T,A,\theta,\eps)=\left\{
\begin{array}{l}
\dfrac{1}{\eps\omega_\eps}\left(\overline{\theta}_2(T_1(A,\theta,\eps)+\eps T,T_1(A,\theta,\eps),\overline{A}_1(T_1(A,\theta,\eps),A,\theta,\eps),\right.\\
\qquad \left.\dfrac{\pi}{2}(1-\eps\omega_\eps),\eps)-\dfrac{\pi}{2}(1-\eps\omega_\eps)\right)-\pi, \qquad {\rm if\ \ }\eps>0,\\
\dfrac{1}{\omega_0}T-\pi, \qquad {\rm if\ \ }\eps=0.
\end{array}
\right.
\end{equation*}
Let us convince ourself that the function $F_2$ verifies the
assumptions of the implicit function theorem \cite[Ch.~X, \S~2,
Theorems~1 and 2]{kolm} at the point $
(T,A,\theta,\eps)=(\omega_0\pi,A_0,\theta_0,\eps). $ Since\\
\centerline{$\dfrac{\pi}{2}(1-\eps\omega_\eps)=\overline{\theta}_2\left(T_1(A,\theta,\eps),T_1(A,\theta,\eps),\overline{A}_1(T_1(A,\theta,\eps),A,\theta,\eps),\dfrac{\pi}{2}(1-\eps\omega_\eps),\eps\right)$}
then
\begin{eqnarray*}
\lim_{\eps \rightarrow 0} F_2(T,A,\theta,\eps)&=&\lim_{\eps \rightarrow 0} \dfrac{1}{\eps\omega_\eps}(\overline{\theta}_2)'_{(1)}(T_1(A,\theta,\eps)+\lambda (A,\theta,\eps)\eps T,T_1(A,\theta,\eps),\\
&&\left.\overline{A}_1(T_1(A,\theta,\eps),A,\theta,\eps),\dfrac{\pi}{2}(1-\eps\omega_\eps),\eps\right)\eps
T-\pi=\dfrac{1}{\omega_0}T-\pi,
\end{eqnarray*}
that is $F_2$ is continuous at $\eps=0$. Here $\lambda (A,\theta,\eps) \in [0,1].$ Furthermore, we have
\begin{equation*}
(F_2)'_T\left(\dfrac{\pi}{2}-\theta_0+\pi,A_0,\theta_0,0\right)=\dfrac{1}{\omega_0}\neq
0.
\end{equation*}
Therefore, the implicit function theorem \cite[Ch.~X, \S~2,
Theorems~1 and 2]{kolm} allows as to conclude that there exists
$\widetilde{T}_2 \in C^1(\mathbb{R} \times \mathbb{R} \times
[0,1],\mathbb{R})$ such that
\begin{eqnarray*}
&&\dfrac{1}{\eps\omega_\eps}\left(\overline{\theta}_2(T_1(A,\theta,\eps)+\eps\widetilde{T}_2(A,\theta,\eps),T_1(A,\theta,\eps),A_1(T_1(A,\theta,\eps),A,\theta,\eps),\right.\\
&&\left.\dfrac{\pi}{2}(1-\eps\omega_\eps),\eps)-\dfrac{\pi}{2}(1-\eps\omega_\eps)\right)=\pi,
\qquad |A-A_0|<\delta,\ |\theta-\theta|<\delta,\ \eps \in
[0,\delta).
\end{eqnarray*}
Since $\theta_0 \in (0,\pi)$, then $\delta>0$ can be diminished in
such a way that
\begin{eqnarray*}
&&\overline{\theta}_3(\pi,T_1(A,\theta,\eps)+\eps\widetilde{T}_2(A,\theta,\eps),\overline{A}_2(T_1(A,\theta,\eps)+\eps \widetilde{T}_2(A,\theta,\eps),\\
&&\left.T_1(A,\theta,\eps),\overline{A}_1(T_1(A,\theta,\eps),A,\theta,\eps),\dfrac{\pi}{2}(1-\eps\omega_\eps),\eps),\dfrac{\pi}{2}(1+\eps\omega_\eps),\eps\right)<\\
&&<\dfrac{\pi}{2}-\theta+\pi+\pi, \qquad \mbox{for}\ \ \eps>0 \ \
\mbox{sufficiently small}
\end{eqnarray*}
and any $|A-A_0|<\delta, |\theta-\theta_0|<\delta$. Therefore, the
solution of (\ref{S3}) with initial conditions under consideration
does not meet the line of discontinuity during time
$(T_1(A,\theta,\eps)+\eps\widetilde{T}_2(A,\theta,\eps),\pi]$.

\begin{figure}\label{fig2}
\begin{center}
\includegraphics[scale=0.9]{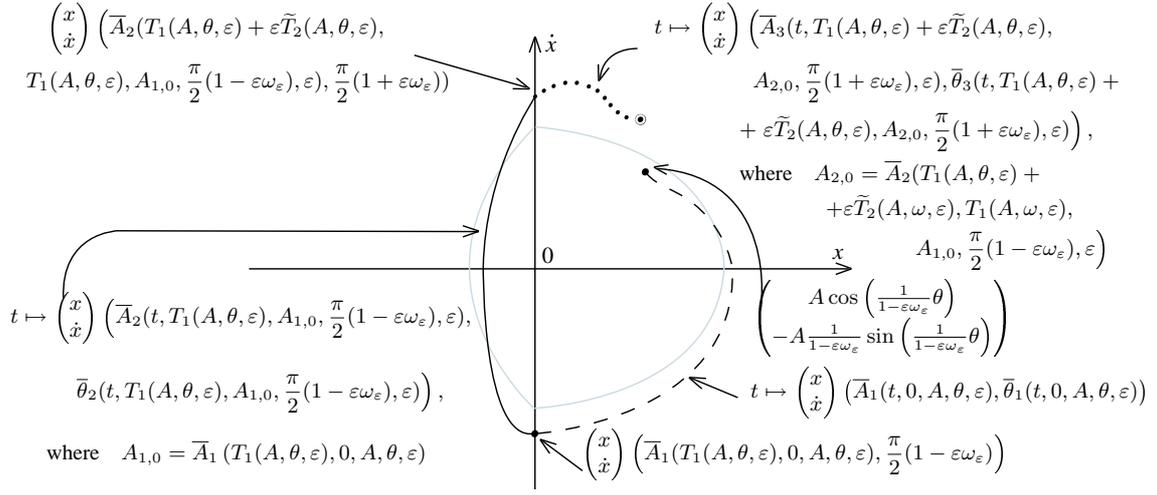}
\footnotesize\caption{Illustration of sewing of the solution of
system (\ref{S1})-(\ref{S1}) by solutions of each its single
component (\ref{S1}), (\ref{S2}) and (\ref{S3}).}
\end{center}
\end{figure}

Summarizing, we can define the solution $t \mapsto
\left(\overline{A}(t,A,\theta,\eps),\overline{\theta}(t,A,\theta,\eps)\right)$
of system (\ref{S1})-(\ref{S3}) as follows (see Fig.~\ref{fig2})
\begin{equation*}
\begin{pmatrix} \overline{A}(t,A,\theta,\eps)\\ \overline{\theta}(t,A,\theta,\eps) \end{pmatrix}=\left\{ \begin{array}{l}
\begin{pmatrix} \overline{A}(t,A,\theta,\eps)\\ \overline{\theta}(t,A,\theta,\eps) \end{pmatrix},\qquad \mbox{if} \quad t \in [0,T_1(A,\theta,\eps)],\\
\begin{pmatrix} \overline{A}_2\left(t,T_1(A,\theta,\eps),\overline{A}(T_1(A,\theta,\eps),A,\theta,\eps),\dfrac{\pi}{2}(1-\eps\omega_\eps),\eps\right) \\ \overline{\theta}_2\left(t,T_1(A,\theta,\eps),\overline{A}(T_1(A,\theta,\eps),A,\theta,\eps),\dfrac{\pi}{2}(1-\eps\omega_\eps),\eps\right)\end{pmatrix},\\
\qquad \qquad \mbox{if} \quad t \in [T_1(A,\theta,\eps),T_1(A,\theta,\eps)+\eps \widetilde{T}_2(A,\theta,\eps)],\\
\begin{pmatrix} \overline{A}_3(t,T_1(A,\theta,\eps)+\eps\widetilde{T}_2(A,\theta,\eps),\overline{A}(T_1(A,\theta,\eps)+\eps \widetilde{T}_2(A,\theta,\eps),A,\theta,\eps),\\ \qquad \qquad \left.\dfrac{\pi}{2}(1+\eps\omega_\eps),\eps\right)\\
\overline{\theta}_3(t,T_1(A,\theta,\eps)+\eps\overline{T}_2(A,\theta,\eps),\overline{A}(T_1(A,\theta,\eps)+\eps
\widetilde{T}_2(A,\theta,\eps),A,\theta,\eps),\\ \qquad \qquad
\left.\dfrac{\pi}{2}(1+\eps\omega_\eps),\eps\right)
\end{pmatrix},\\
\qquad \qquad \mbox{if} \quad t \in [T_1(A,\theta,\eps)+\eps
\widetilde{T}_2(A,\theta,\eps),\pi].
\end{array}
\right.
\end{equation*}

{\bf Step 2.} At this step we show that fixed points of the Poincar\'e map (\ref{PO}) can be studied by means of the function $\overline{P}$ introduced in the formulation of the theorem. To this end we decompose $P_\eps$ as
\begin{equation*}
P_\eps(a,\theta)=\begin{pmatrix} A\\ \theta \end{pmatrix} +\eps(\overline{P}_{\eps,1}(A,\theta)+\overline{P}_{\eps,2}(A,\theta)+\overline{P}_{\eps,3}(A,\theta)),
\end{equation*}
where
\begin{eqnarray*}
&&\overline{P}_{\eps,1}(A,\theta)=\int\limits_0^{T_1(A,\theta,\eps)} G_1(\tau,\overline{A}(\tau,A,\theta,\eps),\overline{\theta}(\tau,A,\theta,\eps),\eps)d\tau,\\
&&\overline{P}_{\eps,2}(A,\theta)=\int\limits_{T_1(A,\theta,\eps)}^{T_1(A,\theta,\eps)+\eps\widetilde{T}_2(A,\theta,\eps)} G_2(\tau,\overline{A}(\tau,A,\theta,\eps),\overline{\theta}(\tau,A,\theta,\eps),\eps)d\tau,\\
&&\overline{P}_{\eps,3}(A,\theta)=\int\limits_{T_1(A,\theta,\eps)+\eps
\widetilde{T}_2(A,\theta,\eps)}^\pi
G_3(\tau,\overline{A}(\tau,A,\theta,\eps),\overline{\theta}(\tau,A,\theta,\eps),\eps)d\tau.
\end{eqnarray*}
Since $\sin, \cos$ and $g$ are bounded on any bounded set then
from system (\ref{S1})-(\ref{S3}) we have that
\begin{equation*}
\begin{pmatrix} \overline{A}(t,A,\theta,\eps)\\ \overline{\theta}(t,A,\theta,\eps) \end{pmatrix} \to  \begin{pmatrix} A\\ t+\theta \end{pmatrix} \quad \mbox{as} \quad \eps \to 0
\end{equation*}
uniformly with respect to $t \in [0,\pi], \ |A-A_0|<\delta, \ |
\theta-\theta_0|<\delta.$ This gives us immediately that
\begin{equation}\label{COP}
\begin{array}{l}
\overline{P}_{\eps,1}(A,\theta) \to \int\limits_0^{T_1(A,\theta,0)} G_1(\tau,A,\tau+\theta,0)d\tau,\\
\overline{P}_{\eps,3}(A,\theta) \to
\int\limits_{T_1(A,\theta,0)}^\pi
G_3(\tau,A,\tau+\theta,0)d\tau,\quad \mbox{as} \quad \eps \to 0,
\end{array}
\end{equation}
uniformly with respect to $|A-A_0|<\delta,$
$|\theta-\theta_0|<\delta.$ Since we proved that $T_1$ and
$\widetilde{T}_2$ are continuously differentiable, then
(\ref{COP}) implies that
\begin{eqnarray*}
&&(\overline{P}_{\eps,_1})'(A,\theta) \to (P_{0,1})'(A,\theta),\\
&&(\overline{P}_{\eps,_3})'(A,\theta) \to (P_{0,3})'(A,\theta),
\quad \mbox{as} \quad \eps \to 0,
\end{eqnarray*}
uniformly with respect to $|A-A_0|<\delta,$
$|\theta-\theta_0|<\delta$.

Let us now study the behavior of $\overline{P}_{\eps,2}$ and $(\overline{P}_{\eps,2})'$ as $\eps \to 0.$ We have
\begin{eqnarray*}
\overline{P}_{\eps,2}(A,\theta)&=&-(1-\eps\omega_\eps) \int\limits_{T_1(A,\theta,\eps)}^{T_1(A,\theta,\eps)+\eps\widetilde{T}_2(a,\theta,\eps)} \begin{pmatrix} \dfrac{1}{\eps} \sin \left(\dfrac{1}{\eps\omega_\eps}\left(\overline{\theta}(\tau,A,\theta,\eps)-\dfrac{\pi}{2}(1-\eps\omega_\eps)\right)+\dfrac{\pi}{2}\right)\\\dfrac{1}{\overline{A}(\tau,A,\theta,\eps)}(1-\eps\omega_\eps)\omega_\eps \cos\left(\dfrac{1}{\eps\omega_\eps}(\overline{\theta}(\tau,A,\theta,\eps)-\right.\\
\qquad \qquad\left.\left.-\dfrac{\pi}{2}(1-\eps\omega_\eps)\right)+\dfrac{\pi}{2}\right) \end{pmatrix} \cdot\\
&&\cdot g\left(\tau,\eps \overline{A}(\tau,A, \theta,\eps)\dfrac{\omega_\eps}{1-\eps\omega_\eps}\cos \left(\dfrac{1}{\eps\omega_\eps}\left(\overline{\theta}(\tau,A,\theta,\eps)-\dfrac{\pi}{2}(1-\eps\omega_\eps)\right)+\dfrac{\pi}{2}\right)\right.,\\
&&\left.-\overline{A}(\tau,A,\theta,\eps)\dfrac{1}{1-\eps\omega_\eps}\sin\left(\dfrac{1}{\eps\omega_\eps}\left(\overline{\theta}(\tau,A,\theta,\eps)-\dfrac{\pi}{2}(1-\eps\omega_\eps)\right)+\dfrac{\pi}{2}\right),\eps\right)d\tau.
\end{eqnarray*}
Scaling the time in the integral as $\tau=T_1(A,\theta,\eps)+\eps\omega_\eps s$ we get
\begin{eqnarray*}
&&\overline{P}_{\eps,2}(A,\theta)=-\omega_\eps(1-\eps\omega_\eps) \int\limits_0^{T_2(A,\theta.\eps)\backslash \omega_\eps} \left(\begin{array}{l} \sin \left(\dfrac{1}{\eps\omega_\eps}\left(\overline{\theta}(T_1(A,\theta,\eps)+\eps s \omega_\eps,A,\theta,\eps)-\right.\right.\\
\left.\left.\quad-\dfrac{\pi}{2}(1-\eps\omega_\eps)\right)+\dfrac{\pi}{2}\right)\\
\eps \dfrac{1}{\overline{A}(T_1(A,\theta,\eps)+\eps\omega_\eps
s,A,\omega,\eps)}(1-\eps\omega_\eps)\omega_\eps \cdot\\
\quad\cdot\cos
\left(\dfrac{1}{\eps\omega_\eps}\left(\overline{\theta}(T_1(A,\theta,\eps)+\eps\omega_\eps
s,A,\theta,\eps)-\right.\right.\\
\qquad \left.\left.-\dfrac{\pi}{2}(1-\eps\omega_\eps)\right)+\dfrac{\pi}{2}\right) \end{array}\right) \cdot \\
&&\cdot g\left(T_1(A,\theta,\eps)+\eps\omega_\eps s,\eps \overline{A}(T_1(A,\theta,\eps)+\eps\omega_\eps s,A,\theta,\eps)\dfrac{\omega_\eps}{1-\eps\omega_\eps} \cos \left(\dfrac{1}{\eps\omega_\eps}(\overline{\theta}(T_1(A,\theta,\eps)+\right.\right.\\
&&\left.\left.+\eps\omega_\eps s,A,\theta,\eps)-\dfrac{\pi}{2}(1-\eps\omega_\eps)\right),-\overline{A}(T_1(A,\theta,\eps)+\eps\omega_\eps s,A,\theta,\eps)\dfrac{1}{1-\eps\omega_\eps}\cdot\right.\\
&&\left.\cdot\sin\left(\dfrac{1}{\eps\omega_\eps}\left(\overline{\theta}(T_1(A,\theta,\eps)+\eps\omega_\eps
s,A,\theta,\eps)-\dfrac{\pi}{2}(1-\eps\omega_\eps)\right)+\dfrac{\pi}{2}\right),\eps\right)ds.
\end{eqnarray*}

\noindent Put
\begin{equation*}
K_\eps(A,\theta)=
\dfrac{1}{\eps}\left(\overline{\theta}(T_1(A,\theta,\eps)+\eps\omega_\eps
s,A,\theta,\eps)-\overline{\theta}(T_1(A,\theta,\eps),A,\theta,\eps)\right).
\end{equation*}
Since
\begin{eqnarray*}
&&\dfrac{1}{\eps}\left(\overline{\theta}(T_1(A,\theta,\eps)+\eps\omega_\eps s,A,\theta,\eps)-\dfrac{\pi}{2}(1-\eps\omega_\eps)\right)=\\
&&=K_\eps(A,\theta) \to
(\overline{\theta})'_{(1)}(T_1(A,\theta,0),A,\theta,0)\omega_0
s=\omega_0 s, \quad \mbox{as} \ \eps \to 0,
\end{eqnarray*}
then
\begin{equation*}
\overline{P}_{\eps,2}(A,\theta) \to -\omega_0 \int\limits_0^\pi
\begin{pmatrix} \sin\left(s+\dfrac{\pi}{2}\right)\\0 \end{pmatrix}
g\left(\dfrac{\pi}{2}-\theta,0,A\sin\left(s+\dfrac{\pi}{2}\right),0\right)ds,
\quad \mbox{as} \ \ \eps \to 0,
\end{equation*}
uniformly with respect to $|A-A_0|<\delta,|\theta-\theta_0|<\delta.$
Since $(K_\eps)'(A,\theta)$ converges as $\eps \to 0$ uniformly in $|A-A_0|<\delta$ and $|\theta -\theta_0|<\delta$ then $(K_\eps)'(A,\theta) \to (K_0)'(A,\theta)$ as $\eps \to 0.$ Therefore,
\begin{equation*}
(\overline{P}_{\eps,2})'(A,\theta) \to (\overline{P}_{0,2})'(A,\theta) \quad \mbox{as} \quad \eps \to 0
\end{equation*}
uniformly with respect to $|A-A_0|<\delta,|\theta-\theta_0|<\delta.$

Summarizing, we proved, that
\begin{eqnarray*}
&&\dfrac{1}{\eps}\left(P_\eps(A,\theta)-\begin{pmatrix} A\\ \theta \end{pmatrix}\right)= \overline{P}(A,\theta)+\\
&&+(\overline{P}_{\eps,1}-\overline{P}_{0,1}(A,\theta)+\overline{P}_{\eps,2}(A,\theta)-\overline{P}_{0,2}(A,\theta)+\overline{P}_{\eps,3}(A,\theta)-\overline{P}_{0,3}(A,\theta)).
\end{eqnarray*}
Therefore, for any $\eps>0$ sufficiently small, the function
\begin{equation*}
(A,\theta) \mapsto P_\eps(A,\theta)-\begin{pmatrix} a\\ \theta \end{pmatrix}
\end{equation*}
has a unique zero $(A_\eps,\theta_\eps)$ such that
$(A_\eps,\theta_\eps)\to (A_0,\theta_0)$ as $\eps \to 0$ and the
real parts of eigenvalues of $(P_\eps)'(A_\eps,\theta_\eps)-I$ are
negative. This is equivalent to say that the eigenvalues of the
matrix $((P_\eps)'(A_\eps,\theta_\eps))^2$ belong to the interval
$[0,1)$. This implies (see \cite[Lemma~9.2]{kraope}) that
\linebreak $t \mapsto
(\overline{A}(t,A_\eps,\theta_\eps,\eps),\overline{\theta}(t,A_{\eps},\theta_\eps,\eps)-t)$
is an asymptotically stable $\pi$-periodic solution of
(\ref{S1})-(\ref{S3}). To see the latter one will probably wish to
make the change of variables $\Xi(t)=\theta(t)-t$ in
(\ref{S1})-(\ref{S3}). Since the change of variables
(\ref{CH1})-(\ref{CH3}) is $\pi$-periodic, than given by these
formulas corresponding solution $(x_\eps,\dot{x}_\eps)$ of
(\ref{EQ1})-(\ref{EQ2}) is also $\pi$-periodic and asymptotically
stable. Uniqueness of $x_\eps$ follows from uniqueness of
$(A_\eps,\theta_\eps)$ and the fact that the change of variables
(\ref{CH1})-(\ref{CH3}) is one-to-one.

To finish the proof it remains to observe that for any $t \in
\left[0,\dfrac{\pi}{2}-\theta_0\right)$ we have
\begin{equation*}
\begin{pmatrix} \overline{A}(t,A_\eps,\theta_\eps,\eps)\cos \left(\dfrac{1}{1-\eps\omega_\eps}\overline{\theta}(t,A_\eps,\theta_\eps,\eps)\right)\\  -A(t,A_\eps,\theta_\eps,\eps)\dfrac{1}{1-\eps\omega_\eps}\sin \left(\dfrac{1}{1-\eps\omega_\eps}\overline{\theta}(t,A_\eps,\theta_\eps,\eps)\right) \end{pmatrix}  \to \begin{pmatrix} A_0\cos (t+\theta_0)\\ -A_0 \sin(t+\theta_0) \end{pmatrix} \quad \mbox{as} \quad \eps \to 0
\end{equation*}
and that for any $t \in \left(\dfrac{\pi}{2}-\theta_0,\pi\right]$
we have
\begin{equation*}
\begin{pmatrix} \overline{A}(t,A_\eps,\theta_\eps,\eps)\cos \left(\dfrac{1}{1-\eps\omega_\eps}\left(\overline{\theta}(t,A_\eps,\theta_\eps,\eps)-\dfrac{\pi}{2}(1+\eps\omega_\eps)\right)+\dfrac{\pi}{2}+\pi\right)\\  -A(t,A_\eps,\theta_\eps,\eps)\sin \left(\dfrac{1}{1-\eps\omega_\eps}\left(\overline{\theta}(t,A_\eps,\theta_\eps,\eps)-\dfrac{\pi}{2}(1+\eps\omega_\eps)\right)+\dfrac{\pi}{2}+\pi\right) \end{pmatrix}  \to
\end{equation*}
\begin{equation*}
\to\begin{pmatrix} A_0\cos (t+\theta_0+\pi)\\ -A_0
\sin(t+\theta_0+\pi) \end{pmatrix} \quad \mbox{as} \quad \eps \to
0,
\end{equation*}
that is $x_\eps$ converges to the solution of (\ref{np}) with the
initial condition $x_\eps(0)=(A_0\cos \theta_0,-A_0\sin \theta_0)$
as $\eps \to 0$ pointwise on $[0,\pi]\backslash\{t_0\}$. The proof
is complete.

\qed

\section {An application}
In this section we apply the result of section 2 to an
impact oscillator shown in Fig. 3.

\parbox[c]{60mm}{
\includegraphics[scale=0.5]{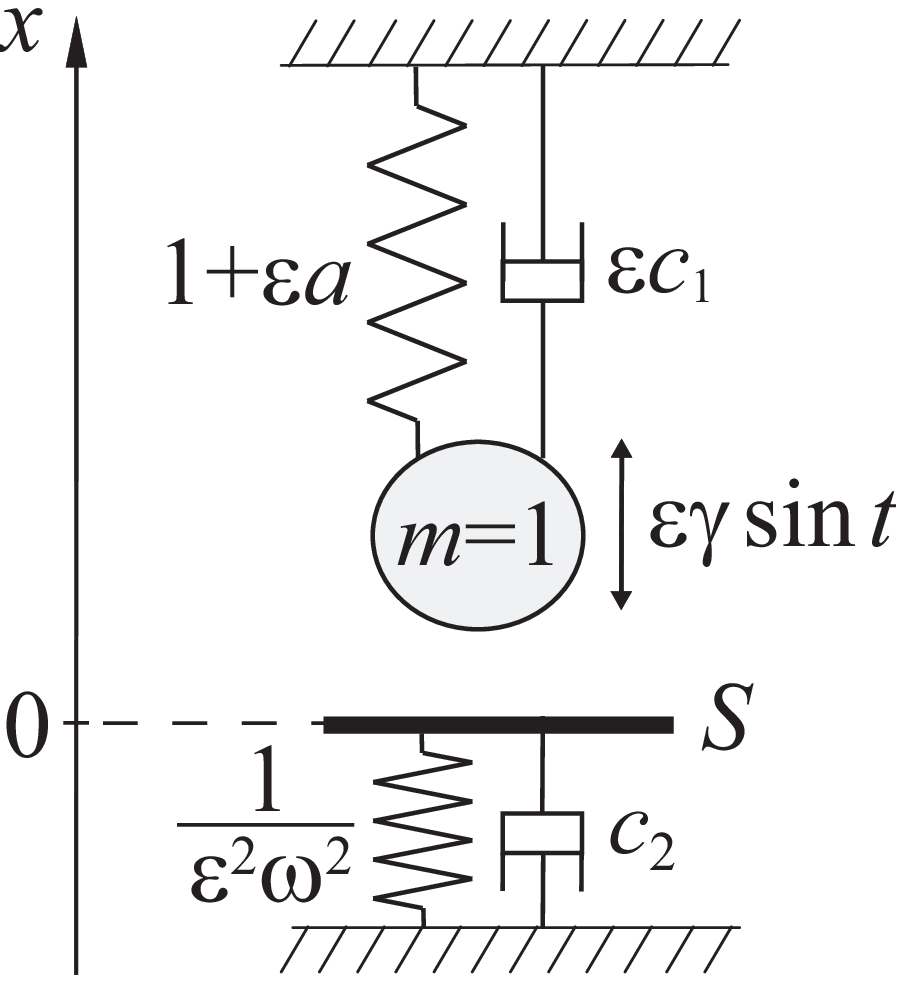}
}\parbox[c]{85mm}{\vskip-1cm \footnotesize Fig. 3: Model of a
preloaded ball (body of mass $m=1$) bouncing against a nearly
elastic surface $S$ which is represented by a spring of stiffness
$\dfrac{1}{\eps^2\omega^2}$ with $\eps>0$ small. The rest
coordinate of the mass is assumed to coincide with the origin of
the  $x$-axis on the one hand and with the surface $S$ at rest on the other
hand. We also suppose that the viscous friction coefficient equals
 $\eps c_1$ outside the contact with $S$ and that it takes the
value $c_2+\eps c_1$ during the contact.}

 \noindent A body of mass $m=1$ is bouncing against a nearly elastic surface $S$
(large stiffness $1/  \eps^2 \omega^2$).
Assuming in addition that the body
is subjected to Rayleigh excitation, viscous friction and forcing, the equation of motions can
be written as follows
\begin{equation}\label{new313}
\begin{array}{ll}
\ddot{x}+x=-\eps a x -\eps c_1 \dot{x} + \eps\mu_1\dot x (1-\dot x^2)+\eps \gamma \sin t, & \mbox{if} \quad x \geqslant 0,\\
\ddot{x}+\dfrac{1}{\eps^2\omega^2}x=-(c_2+\eps
c_1)\dot{x}+(\mu_2+\eps\mu_1)\dot x(1-\dot x^2)+\eps \gamma \sin
t, & \mbox{if} \quad x<0.
\end{array}
\end{equation}

\noindent System (\ref{new313}) is of the form (\ref{ps}) with limiting system. There is
no obvious reason why the Rayleigh excitation should be $O(1)$ during impact, but as
this does not complicate the analysis, we admit this possibility.
The averaging function $\overline{P}$ takes now the longer
form
\begin{eqnarray*}
\overline{P}(A,\theta)&=&-\int_0^{\frac{\pi}{2}-\theta} \begin{pmatrix} \sin(\tau+\theta)\\ \frac{1}{A}\cos(\tau+\theta) \end{pmatrix} (-aA\cos(\tau+\theta)+c_1A\sin(\tau+\theta)+\\
&&-\mu_1 A\sin (\tau+\theta)(1-A^2\sin^2 (\tau+\theta))+\gamma \sin\tau-2\omega A\cos(\tau+\theta))d\tau-\\
&&-\int_{\frac{\pi}{2}-\theta}^{\pi} \begin{pmatrix} \sin(\tau+\theta+\pi)\\ \frac{1}{A}\cos(\tau+\theta+\pi) \end{pmatrix} (-aA\cos(\tau+\theta+\pi)+c_1A\sin(\tau+\theta+\pi)+\\
&&-\mu_1 A\sin(\tau+\theta+\pi)(1-A^2\sin^2(\tau+\theta+\pi))+\gamma \sin\tau-2\omega A\cos(\tau+\theta+\pi))d\tau-\\
&&-\omega \int_0^\pi \begin{pmatrix} \sin(\tau+\frac{\pi}{2})\\0 \end{pmatrix} \left(c_2 A \sin \left(\tau+\frac{\pi}{2}\right)-\mu_2A\sin\left(\tau+\frac{\pi}{2}\right)\left(1-A^2\sin^2\left(\tau+\frac{\pi}{2}\right)\right)\right)d\tau=\\
&&=\begin{pmatrix} \gamma \theta \cos \theta -\frac{\pi}{2}A(c_1+c_2 \omega)+ \frac{\pi}{2}(\mu_1+\mu_2\omega)A\left(1-\frac{3}{4}A^2\right)\\ -\frac{1}{A} \gamma \cos\theta-\frac{1}{A}\gamma\theta\sin\theta+\frac{\pi}{2}(a+2\omega) \end{pmatrix}.
\end{eqnarray*}
To formulate our result we need some preliminary notations. First, we introduce the function $M:\mathbb{R}\to\mathbb{R}$ as follows
\begin{eqnarray*}
  M(\theta)&=&-\theta\cos\theta+(c_1+c_2\omega)\frac{\cos\theta+\theta\sin\theta}{a+2\omega}+\\
  &&-(\mu_1+\mu_2\omega)\left(\frac{\cos\theta+\theta\sin\theta}{a+2\omega}-\frac{3\gamma^2}{\pi^2}\left(\frac{\cos\theta+\theta\sin\theta}{a+2\omega}\right)^3\right).
\end{eqnarray*}
\begin{proposition}\label{B}
Let $M(\theta_0)=0$ for some $\theta_0 \in (0,\frac{\pi}{2})$ and $A_0$ is defined as
\begin{equation}\label{A0}
  A_0=\frac{\gamma\cos\theta_0+\gamma\theta_0\sin\theta_0}{\frac{\pi}{2}(a+2\omega)}.
\end{equation}
If
\begin{equation}\label{Mprime}
M'(\theta_0)>0
\end{equation}
and
\begin{equation*}
-(c_1+c_2\omega)+(\mu_1+\mu_2\omega)(1-3{A_0}^2)<0
\end{equation*}
then, for any $\eps>0$ sufficiently small, equation (\ref{new313})
has exactly one $\pi$-periodic solution $x_\eps$ such that
\begin{equation*}
  (x_\eps(0),\dot x_\eps(0)) \to (A_0\cos\theta_0,-A_0\sin\theta_0) \ \mbox{as} \ \eps \to 0.
\end{equation*}
The solution $x_\eps$ is asymptotically stable.
\end{proposition}
The proof of proposition (\ref{B}) relies on the following lemma.
\begin{lemma}\label{C}
Consider $2\times2$ real matrix $D$. If $Sp \ D<0$ and
$\det\|D\|>0$ then the eigenvalues of $D$ have negative real
parts.
\end{lemma}
The statement of the lemma follows from the direct computation of the eigenvalues of $D$ according to the standard formula for roots of quadratic equations.\\

\noindent{\bf Proof of proposition \ref{B}.} Direct computation
shows that $(A_0,\theta_0)$ is a zero of $\overline{P}$. To prove
the proposition it remains to show that
\begin{itemize}
\item[1)] $Sp \ \overline P'(A_0,\theta_0)<0$.
\item[2)] $\det \|\overline P'(A_0,\theta_0)\|>0$.
\end{itemize}
But these two relations follow from the formulas
\begin{itemize}
\item[1)] $Sp \ \overline P'(A_0,\theta_0)=-\pi(c_1+c_2\omega)+\pi(\mu_1+\mu_2\omega)(1-3A_0^2) $.
\item[2)] $\dfrac{2A_0}{\pi\gamma(a+2\omega)}\det \|\overline P'(A_0,\theta_0)\|=M'(\theta_0)$,
\end{itemize}
which are straightforward.

\qed

 Our next proposition \ref{A} shows that proposition \ref{B}
is not vacuous, namely we give sufficient conditions ensuring that
(\ref{Mprime}) is satisfied. Before proceeding to the formulation of
proposition \ref{A} we need to introduce some notations and
properties. First we observe, that the
\begin{equation*}
  \frac{M(\theta)}{\theta\cos\theta}=-1+K(\theta),
\end{equation*}
where
\begin{equation*}
  K(\theta)=\left(\frac{c_1+c_2\omega}{a+2\omega}-\frac{\mu_1+\mu_2\omega}{a+2\omega}\left(1-\frac{3\gamma^2}{\pi^2}\left(\frac{\cos\theta+\theta\sin\theta}{a+2\omega}\right)^2\right)\right)\left(\frac{1}{\theta}+ tg \ \theta\right).
\end{equation*}
Observe, that there exists $\theta_* \in \left(0,\frac{\pi}{2}\right)$ such that
\begin{equation}\label{KK}
  K'(\theta) \ne 0 \ \mbox{for all} \ \theta \in \left(\theta_*,\frac{\pi}{2}\right).
\end{equation}
In fact, if $K'(\theta_n)=0$ for some sequence $\theta_n \uparrow
\frac{\pi}{2}$ as $n \to \infty$, then
\begin{equation*}
  M'(\theta_n)=(\cos\theta_n-\theta_n\sin\theta_n)(-1+K(\theta_n)) \to \infty \ \mbox{as} \ n \to \infty,
\end{equation*}
that contradicts the boundedness of the derivative of $M$.
\begin{proposition}\label{A}
  Let $\theta_*\in\left(0,\frac{\pi}{2}\right)$ be such a number that (\ref{KK}) holds true. Assume that
\begin{equation*}
  K(\theta_*)<1
\end{equation*}
and denote by $\theta_0 \in (\theta_*,\frac{\pi}{2})$ the unique point satisfying
\begin{equation}\label{sati}
  -1+K(\theta_0)=0.
\end{equation}
Then $M(\theta_0)=0$ and $M'(\theta_0)<0$. Consequently, for any
$\eps>0$ sufficiently small, equation (\ref{new313}) has
exactly one $\pi$-periodic solution $x_\eps$ such that
\begin{equation*}
  (x_\eps(0),\dot x_\eps(0)) \to (A_0\cos\theta_0,-A_0\sin\theta_0) \ \mbox{as} \ \eps \to 0,
\end{equation*}
where
\begin{equation*}
  A_0=\frac{\gamma\cos\theta_0+\gamma\theta_0\sin\theta_0}{\frac{\pi}{2}(a+2\omega)}.
\end{equation*}
The solution $x_\eps$ is asymptotically stable.
\end{proposition}
Note that $\theta_0 \in \left(\theta_*,\frac{\pi}{2}\right)$
satisfying (\ref{sati}) always exists and is unique since
$-1+K(\theta_*)<0,\ K'(\theta)\ne 0, \ \mbox{for} \ \theta \in
\left(\theta_*,\frac{\pi}{2}\right),$ and
\begin{equation}\label{KKK}
  \lim\limits_{t\uparrow \frac{\pi}{2}} K(\theta)=+\infty.
\end{equation}

\noindent{\bf Proof of proposition \ref{A}.} First, observe that
$M(\theta_0)=\theta_0\cos\theta_0(-1+K(\theta_0))=0$. Second, we
have
\begin{equation*}
  M'(\theta)=(\cos\theta-\theta\sin\theta)L(\theta)+\theta\cos\theta K'(\theta)
\end{equation*}
and so $M'(\theta_0)=K'(\theta_0)$. But properties (\ref{KK}) and
(\ref{KKK}) imply that $K'(\theta_0)>0$ and so the proof is
complete.

\qed

Let us finally formulate our result in the simpler setting when the
Rayleigh excitation is switched off, that is $\mu_1=\mu_2=0$. We
have that
\begin{equation*}
  K(\theta)=\frac{c_1+c_2\omega}{a+2\omega}\left(\frac{1}{\theta}+ tg \
  (\theta)\right),
\end{equation*}
in particular, there exists an unique $\theta_*
\in\left(0,\frac{\pi}{2}\right)$ such that $sign \
K'(\theta_*)=sign \ (\theta-\cos\theta)=0$.
\begin{proposition}\label{D}
Let $\theta_* \in\left(0,\frac{\pi}{2}\right)$ be the unique point such that $\theta_*-\cos\theta_*=0$. Assume that
\begin{equation*}
  \frac{c_1+c_2\omega}{a+2\omega}\left(\frac{1}{\theta_*}+tg \ \theta_*\right)<1.
\end{equation*}
Denote by $\theta_0 \in\left( \theta_*,\frac{\pi}{2} \right)$ the unique point such that
\begin{equation*}
  -1+\frac{c_1+c_2\omega}{a+2\omega}\left(\frac{1}{\theta_0}+tg \ \theta_0\right)=0.
\end{equation*}
Then the conclusion of proposition \ref{A} holds true.
\end{proposition}

\section{Discussion}
\begin{itemize}
\item It is remarkable that the analysis of system (\ref{ps}) can be handled by the introduction
of a Poincar\'e map and the use of the implicit function theorem although the limiting system
for $\varepsilon \rightarrow 0$ is non-smooth.
\item
At the same time we have formulated an unusual type of singular perturbation problem.
Putting $\varepsilon =0$, we have non-smooth impact, for $\varepsilon > 0$ we have fast
motion in a neighborhood of the subset $x=0$. For $x>0$ slow motion takes place but
 this is not described by standard slow manifold theory, see
\cite{V05}. Still, the dynamics for $x>0$ can be considered as taking place in an explicitly
formulated slow manifold.\\
On the other hand, the solutions for
$x <0$ have as slow manifold the boundary $x=0$. This does not satisfy the necessary hyperbolicity
condition, but
the solutions for $x >0$ are forced to the manifold $x=0$ and, after passing by a fast transition
through the domain $x <0$ they are forced again to leave $x=0$. We note also that sliding along the
slow manifold, as happens for instance in dry friction problems, is not possible. This simplifies
the bifurcational behavior.
\item
Regarding the averaging result obtained in this paper, we draw attention to the papers \cite{P79}-
 \cite{MS85} and further references there. In  \cite{P79} the framework of differential inclusions
is used, in \cite{MS85} explicit estimates of the vector field and the solutions are given in the case
of impulsive forces. Our approach economically avoids the estimate of general solution behavior
as we aim at the more restricted result of obtaining periodic solutions.
\end{itemize}


\begin{thebibliography}{99}

\bibitem{BK01} V.I. Babitsky, V.L. Krupenin, ``Vibration of strongly
nonlinear discontinous systems," Springer, Berlin etc., (2001).

\bibitem{blm} A. Buica, J. Llibre, O. Makarenkov,
Asymptotic stability of periodic solutions for nonsmooth
differential equations with application to the nonsmooth van der
Pol oscillator, SIAM J. Math. Anal. 40 (2009), no. 6, 2478--2495.

\bibitem{lazer} J. Glover, A.C. Lazer, P.J. McKenna,
Existence and stability of large scale nonlinear oscillations in
suspension bridges. Z. Angew. Math. Phys. 40 (1989), no. 2,
172--200.

\bibitem{ivanov} A. Ivanov, Bifurcations in impact systems, Chaos, Solitons \& Fractals Vol. 7, No. 10. pp. 1615-1634,
1996.

\bibitem{kolm} A. N. Kolmogorov, S. V. Fom{\rm i}n, ``Elements of the
theory of functions and functional analysis," Fourth edition,
revised. Izdat. ``Nauka'', Moscow, 1976  (in Russian); transl. 1st ed.
Dover Publ., New York (1996).

\bibitem{kraope}
M. A. Krasnosel'skii, ``The operator of translation along the
trajectories of differential equations," Translations of
Mathematical Monographs, 19. Translated from the Russian by
Scripta Technica, American Mathematical Society, Providence, R.I.
(1968).

\bibitem{K08} P.O.K. Krehl, ``Details of history of shock waves, explosions and impact:
a chronological and biographical reference," Springer, Berlin etc., (2008).

\bibitem{mawhin} J. Mawhin,
Resonance and nonlinearity: a survey. Ukrainian Math. J. 59
(2007), no. 2, 197--214.

\bibitem{MS85} Yu.A. Mitropolsky, A.M. Samoilenko, ``Forced oscillations of systems with
impulsive force'', Int. J. Non-Linear Mechanics 20, pp. 419-426 (1985)

\bibitem{P79} V.A. Plotnikov, ``The averaging method for differential inclusions and its
application to optimal-control problems'', Differential Equations 15, pp. 1427-1433 (1979)

\bibitem{V05} F. Verhulst, ``Methods and applications of singular perturbations,
boundary layers and multiple timescale dynamics," Springer, New York etc., (2005).

\bibitem{yagasaki} K. Yagasaki, Nonlinear dynamics of vibrating microcantilevers in tapping-mode atomic force
microscopy, Physical Review B 70, 245419 (2004).

\end{thebibliography}
\end{document}